\documentclass[a4paper,12pt]{article}
\usepackage[latin1]{inputenc}
\usepackage[T1]{fontenc}
\usepackage[francais, english]{babel}
\usepackage{amssymb,amsmath}
\usepackage{amsthm}
\usepackage{mathrsfs}
\usepackage{geometry}

\title{\textbf{ On orbit Reflexive tuple of operators And Weak Orbit Reflexivity  } }\vspace{1cm}
\author{\hspace{0.5cm} $^{1}$Abdelaziz  Tajmouati \hspace{0.5cm} $^{2}$ Youness   Zahouan }

\date{}

\begin{document}
\maketitle
\begin{center}
$^{1,2}$  University Sidi Mohamed Ben Abdellah

 Faculty of Sciences Dhar Al Mahraz, Laboratory of Mathematical Analysis and Applications, Fez, Morocco.\\
Email  : abdelaziztajmouati@yahoo.fr \hspace{0.3cm} zahouanyouness1@gmail.com \\
\end{center}

\textbf{ Abstract:}
In this paper we give a various conditions for which the tuple  $\mathcal{T} = (T_{1} , T_{2} , ... , T_{n})$ of commutative bounded linear operators on an infinite dimensional ( real , complex ) Banach space X  is orbit reflexive. After we  introduce the notion of weak orbit reflexive operator and we show some results.\\


\begin{center}
\LARGE\textbf{Introduction}
\end{center}

By an n-tuple of operators we mean a finite sequence of length $n$ of commuting  bounded linear operators on a Banach space X. The orbit of a points $x \in X$ under an operator $T \in B(X)$  is the sequence,

\begin{center}
 $\left\{T^{n}x : n = 0 , 1 , ... \right\}$
\end{center}

The notion of orbit-reflexive operators on a Hilbert space was introduced and
studied in [6]. While the reflexivity of operators is connected to the invariant
subspace problem, its natural analogue of orbit reflexivity is in the same way connected
to the problem of existence of closed invariant subsets.\\
Recall that if X is a Banach space and  $B(X)$ the set of all bounded linear operators acting on X,an operator $T \in B(X)$ is said to be reflexive if every operator $A \in B(X)$ belongs to 
 $\left\{p(T) : p \; polynomial\right\}^{-SOT}$ (the closure in the strong operator topology) , whenever $Au \in \left\{p(T)x : p \, polynomial\right\}^{-}$  (the norm closure of the set$\left\{p(T) : p \, polynomial\right\}$ ) for each $x \in X$. \\
Similarly an operator is said to be orbit reflexive if every operator $A \in B(X)$ belongs to $ Orb(T)^{-SOT}$, whenever $Au \in  Orb(T,x)^{-}$ for each $x \in X$(or equivalently, such that every closed subset of X invariant for T is invariant for A ).\\
The orbit reflexivity of many classes of Hilbert space operators was shown in
[6], e.g. for normal operators, contractions, algebraic operators, weighted shifts
and compact operators. Among others, each operator whose spectrum does not intersect the unit circle is orbit reflexive, and a various conditions which insure that an operator acting on a Banach space is orbit reflexive are given in [2]. Other authors are introduce a similar notions (Null-orbit reflexive operators [7], $\mathbb C$-orbit reflexive operators [8], $ \mathbb R-$orbit Reflexive Operators [9]).\\
In this work we  focus on the orbit reflexivity.

\begin{flushleft}
\begin{flushleft}
\section{Operator orbit reflexive.}
\end{flushleft}
\end{flushleft}
  
{\bf{Theorem \,1.1 }}\, [  6 , Theorem 5 ]

Let X be an  Hilbert space and $T \in B(X)$. Then T is orbit reflexive in any of the following cases:

\begin{enumerate} 
\item there is a non empty open subset $U \subset X $such that for each $x \in U$, the
orbit $Orb(T,x)$ is closed.

\item  there is a non empty open subset $U \subset X $ such that for each $x \in U$, $ \left\|T^{n}x\right\| \longrightarrow \infty $ , as $n \longrightarrow \infty$.

\item for each $x \in X$ , $ \left\|T^{n}x\right\| \longrightarrow 0 $.

\item The set $ Orb(T)^{-SOT}$ is countable and strongly compact.

\item  $\sigma(T) \cap \left\{ \lambda \in \mathbb C / \left| \lambda\right| = 1 \right\} = \emptyset $.

\end {enumerate}

{\bf{ Theorem \,1.2 }}\, [2 , Theorem 7 ] T is orbit reflexive if,
\begin{center}
$\sum^{\infty}_{n=1} 1/\left\|T^{n}\right\| < \infty$
\end{center}

or when X is a complex Hilbert space 

\begin{center}
$\sum^{\infty}_{n=1} 1/\left\|T^{n}\right\|^{2}_{} < \infty$
\end{center}

{\bf{ Corollary \,1.1 }}\, [2 , Corollary 8  ]
Let X be an  Hilbert space and $T \in B(X)$. If $ r(T) \neq 1$ then T is orbit reflexive.\\

{\bf{Proof.}}
\begin{itemize}
\item if $r(T) < 1$ then $lim_{n \rightarrow \infty}\left\|T^{n}\right\|=0 $ so , $ \left\|T^{n}x\right\| \longrightarrow 0$ for each $x \in X$ and we can apply Theorem 1.1 (3)
\item if $r(T) > 1$ then $\left\|T^{n}\right\| > n^{2}$ for all $n \in N$ large enough (since otherwise we would have $r(T) =inf_{n \rightarrow \infty} \left\|T^{n}\right\|^{1/n} \leq lim_{n \rightarrow \infty} (n^{2})^{1/n} \leq 1 $) , and we can apply Theorem 1.2
\end{itemize}

{\bf{ Corollary \,1.2 }}\,
Let X be an  Hilbert space and $T \in B(X)$.\\
If $ r(T) \neq 1$, then $T^{*}$ is orbit reflexive.\\

{\bf{Proof.}}\, Follows from, $ r(T) =  r(T^{*})$.\\

{\bf{Corollary \,1.3 }}\,
Let X be a  Banach  space and $T \in B(X)$.If T is self-adjoint and $ r(T) \neq 1$, Then $T^{n}$ is orbit reflexive for all $ n \geq 2$.\\

{\bf{Proof.}}

\begin{itemize}
\item if $n = 2$,
\begin{center}
$\left\|T^{2}\right\| = sup_{\left\|x\right\|=1} \left|\left\langle T^{*}Tx,x\right\rangle\right|=sup_{\left\|x\right\|=1} \left|\left\langle Tx,Tx\right\rangle\right|= sup_{\left\|x\right\|=1}\left\|Tx\right\|^{2}= \left\|T\right\|^{2}$.
\end{center}

And by recursion:
\begin{center}
 $\forall m \in N , \left\|T^{2^{m}}\right\|= \left\|T\right\|^{2^{m}}$
\end{center}

Thus , 
\begin{center}
$r(T)= lim_{m\rightarrow \infty }\left\|T^{m}\right\|^{1/m}= lim_{m\rightarrow \infty }\left\|T^{2^{m}}\right\|^{1/2m}=  \left\|T\right\| $
\end{center}

Then ,
\begin{center}
 $r(T^{2})= \left\|T^{2}\right\|=\left\|T\right\|^{2} = r(T)^{2}.$
\end{center}

Or, $ r(T) \neq 1$ Then $r(T^{2})\neq 1$ and by Corollary 1.1, $T^{2}$ is orbit reflexive.

\item By recursion we have: 
\begin{center}
 $r(T^{n})= r(T)^{n}$ and by Corollary 1.1, $T^{n}$ is orbit reflexive.
\end{center}
\end{itemize}

\begin{flushleft}
\section{n-tuple of operators orbit reflexive.}
\end{flushleft}

{\bf{Definition \,2.1 }}\,[4 ]
Let $  \mathcal{T} = (T_{1} , T_{2} , ... , T_{n})$ be an n-tuple of operators acting on an infinite dimensional Banach space X. We will let, \\

\begin{center}
$\mathcal{F}= \left\{T_{1}^{k_{1}}T_{2}^{k_{2}} ... T_{n}^{k_{n}} : k_{i} \geq 0,
 i = 1,...,n \right\}$
\end{center}

be the semi-group generated by $\mathcal{T}$. For $x \in X$, the orbit of x under the tuple $\mathcal{T}$ is the set

\begin{center}
$Orb(\mathcal{T},x) = \left\{Sx : S \in  \mathcal{F} \right\} $
\end{center}

{\bf{Definition \,2.2 }}\,[12 , Definition 1.2]
The orbit of x under the tuples $\mathcal{T}$ tending to infinity if:
\begin{center}
 $\left\| T_{1}^{k_{1}}T_{2}^{k_{2}} ... T_{n}^{k_{n}}x  \right\|   \rightarrow   \infty $ as $k_{i}   \rightarrow   \infty $ with  $k_{i} \geq 0 $, for all $ i = 1 , ... ,n $
\end{center}

{\bf{Lemma  \,2.1 }} [6]\, 
Let $ A , S_{1} , S_{2} , ...  \in B(X)$ .\\
If the set of vectors $u \in X$ for which
$Au \in \left\{ S_{1}u , S_{2}u , . . .\right\}$ is of second category, then $  A \in \left\{ S_{1} , S_{2} , . . .\right\}$\\

{\bf{Proof.}}\, Let $U$ the above mentioned set of second category, so $U \subset \bigcup^{\infty}_{n=1}Ker(A -S_{n}) $.By Baire category theorem,one of these closed sets has a none-empty interior, so it must be the whole space X, then $A = S_{n}$ for some n.\\

{\bf{Proposition  \,2.1 }}\,
Let $ A , T_{1} , T_{2} , ... , T_{n} \in B(X)$ and $  \mathcal{T} = (T_{1} , T_{2} , ... , T_{n})$ be an n-tuple of commuting operators.\\
If the set of vectors $x \in X$ for which
$Ax \in Orb(\mathcal{T},x)$ is of second category, then $  A \in Orb(\mathcal{T})$\\

{\bf{Proof.}}\\
Follows from Lemma 2.1 \\

{\bf{Theorem \, 2.3}}\,
Let $  T_{1} , T_{2} , ... , T_{n} \in B(X)$ and $  \mathcal{T} = (T_{1} , T_{2} , ... , T_{n})$ be an n-tuple of commuting operators.\\
$  \mathcal{T}$ is orbit reflexive in any of the following cases:

\begin{enumerate} 
\item there is a non empty open subset $U \subset X $such that for each $x \in U$, the
orbit $Orb(\mathcal{T},x)$ is closed,

\item  there is a non empty open subset $U \subset X $ such that for each $x \in U$, $\left\| T_{1}^{k_{1}}T_{2}^{k_{2}} ... T_{n}^{k_{n}}x  \right\|   \rightarrow   \infty $ as $k_{i}   \rightarrow   \infty $ with  $k_{i} \geq 0 $, for all $ i = 1 , ... ,n $,

\item The set  $Orb(\mathcal{T})^{-SOT}$ is countable and strongly compact,
\end {enumerate}

{\bf{Proof.}}\; Throughout the proof , suppose tha t $Ax \in Orb(\mathcal{T},x)^{-}$  for every $x  \in X $. We need  to show that $A \in Orb(\mathcal{T})^{-SOT}$  in each of the cases:
\begin{enumerate} 
\item Follows from Lemma 1.1,

\item  Since $\left\| T_{1}^{k_{1}}T_{2}^{k_{2}} ... T_{n}^{k_{n}}x  \right\|   \rightarrow   \infty $ as $k_{i}   \rightarrow   \infty $ then $Orb(\mathcal{T},x)$ is closed , it  follows from (1) that $\mathcal{T}$ is  orbit reflexive.

\item let $S \in  \mathcal{F} $ and suppose that  $Ax$ $=$ lim $ S^{n_{k}}x$ for some sequence $(n_{k})_{k \geq 0}$ of positive integers.\\
 Since $Orb(\mathcal{T})^{-SOT}$ is strongly compact, the sequence$ (S^{n_{k}})_{k \geq 0}$has a strongly convergent subsequence. Therefore $Ax \in \left\{ Bx  / B \in Orb(\mathcal{T})^{-SOT} \right\}$  for each $x\in X$,
and according to  Lemma 2.1, we have $A \in Orb(\mathcal{T})^{-SOT}$.

\end {enumerate}

{\bf{Example \,1.1 }}\,

Let S be the unilateral forward shift on $\ell^{2}(\mathbb N)$:

\begin{center}
$Se_{n}= e_{n+1}$    ;    $i\geq 1$
\end{center}

where $\left\{e_{n} : n\in \mathbb N  \right\} $ is the standard orthonormal basis for $\ell^{2}( \mathbb N)$.

Given a sequence of positive numbers $(a_{i})_{i\geq 1}$
so that $a_{i} \succ 1 $ for all $i\geq1$ and $a_{i} \rightarrow 1 $ as $i \rightarrow \infty$ and let \\

\begin{center}
$T_{i} = a_{i}S$  ;  $i = 1 , .... , n$
\end{center}

So, of after[12 , Example 1.1] there is a non empty open subset $U \subset X $ such that for each $x \in U$,the orbit of x under the tuple $(T_{1} , T_{2} , ... , T_{3} ) $ tend strongly to infinity then follows from Theorem \, 2.3 (2)  the tuple $(T_{1} , T_{2} , ... , T_{3} ) $ it orbit reflexive.\\

{\bf{Proposition \, 2.2}}\,
Let X an infinite dimensional reflexive Banach space and  $  \mathcal{T} = (T_{1} , T_{2} , ... , T_{n})$ be the n-tuple of operators in $B(X)$bounded below for all $i \geq 1$.\\
If there is $x \in X$  such that  the orbit of x under $T_{i}$ for all $i \geq 1$ tend strongly to infinity then $\mathcal{T}$ is orbit reflexive.\\

{\bf{Proof.}}\\
From Corollary 1.1 [12] the orbit of x under the tuple $\mathcal{T}$ tend strongly to infinity and from 2 of Theorem \, 2.3 $\mathcal{T}$ is orbit reflexive.\\


\section{Operator weakly orbit reflexive.}
By a weak orbit of an operator  $T \in B(X)$ we mean a sequence of the form $\left(\left\langle T^{n}x , x^{*}\right\rangle\right)^{n=0}_{\infty}$ ,where $x \in X$ and $x^{*}_{} \in X^{*}, $ and we write :
\begin{center}
$W-Orb(T,x)=\left(\left\langle T^{n}x , x^{*}\right\rangle\right)^{n=0}_{\infty}$
\end{center}

{\bf{Definition \, 1}}\,
Let X a Banach space ,an operator $T \in B(X)$ is said to be weakly orbit reflexive if every  operator $A \in B(X)$ belongs to $ W-Orb(T)^{-WOT}$ (the closure in the weak operator topology), whenever $\left\langle Ax , x^{*}\right\rangle \in  W-Orb(T,x)^{-W}$(the closure in the weak topology) for each $x \in X $ and $ x^{*} \in X^{*}$.\\

{\bf{Example \, 1}}\,
\begin{itemize}
\item Since the norm topology is strictly stronger than the
weak topology, and so every  operator orbit reflexive  is a weakly orbit reflexive,so  in a Hilbert space a normal operators, contractions, algebraic operators, weighted shifts and compact operators are weak orbit-reflexive.
\item Each Banach space operator have spectral radius different from 1, is weak orbit-reflexive.

\end{itemize}
{\bf{Theorem \, 3.1}}
\begin{enumerate}
\item If T is weakly orbit reflexive and S is invertible operator  then $STS^{-1}$ is also weakly orbit reflexive.
\item If T is  weakly orbit reflexive then $T \oplus T$ is also orbit reflexive.
\end{enumerate}

{\bf{Proof.}}\,

\begin{enumerate}
\item  The image by an operator of a weakly convergent sequence is a sequence weakly 
convergent.Hence $\left\langle SAS^{-1}x , x^{*}\right\rangle \in  W-Orb(STS^{-1},x)^{-W}$  for all $x \in X$ and $ x^{*} \in X^{*}$ is
equivalent to $\left\langle Ax , x^{*}\right\rangle\in  W-Orb(T,x)^{-W}$ for all $x \in X$ and $ X^{\ast} \in X^{*}$, or T is weakly orbit reflexive , then $A \in  W-Orb(T)^{-WOT}$, so $SAS^{-1} \in  W-Orb(STS^{-1})^{-WOT}$.

\item Let$A \in B(X\oplus X)$and suppose that $\left\langle Ax , x^{\ast} \right\rangle \in W-Orb(T\oplus T,x)^{-W}$ for all $ x \in X \oplus X$ and $ x^{\ast} \in X^{\ast}\oplus X^{\ast}$ .

In particular, for any $a \in X $ we have :
\begin{center}
 $\left\langle A(a,0) ,x^{\ast}\right\rangle \in W-Orb(T \oplus T,(a,0))^{-W}  \subset \mathbb C× \left\{0\right\} $
\end{center}
And similarly for the second component, so that both copies of X are A-invariant. Hence A can be written as $B\oplus C$ where $B,C \in B(X)$.\\

Moreover, if $a \in X $then 

$\left\langle (Ba,Ca), x^{\ast}\right\rangle = \left\langle A(x,x), x^{\ast}\right\rangle \in W-Orb(T\oplus T,x)^{-W} \subset \left\{ {(\lambda ,\lambda):
\lambda \in \mathbb C} \right\}$

 So $B =C$. \\
Now since T is weakly orbit reflexive and $Ba \in W-Orb(T,x)^{-W}$
for all $a \in X $, we have $B \in  W-Orb(T)^{-WOT}$, and therefore:\\
\begin{center}
 $A=B\oplus B \in  W-Orb(T \oplus T)^{-WOT}$
\end{center}
\end{enumerate}

{\bf{ Proposition \,3.1 }}\,  T is weakly orbit reflexive if
\begin{center}
$\sum^{\infty}_{n=1} 1/\left\|T^{n}\right\|^{1/2} < \infty$
\end{center}

or when X is a complex Hilbert space 

\begin{center}
$\sum^{\infty}_{n=1} 1/\left\|T^{n}\right\| < \infty$
\end{center}

{\bf{Proof.}}\,

Let$A \in B(X)$ and suppose that $\left\langle Ax , x^{\ast} \right\rangle \in W-Orb(T,x)^{-W}$ for all $ x \in X $ and $ x^{\ast} \in X^{\ast}$ .\\
Suppose that $A \neq T^{n}$  for all $n \in \mathbb N$, otherwise there is nothing to prove.\\

Since $\sum^{\infty}_{n=1} \frac{1}{\left\|T^{n}\right\|^{1/2}} < \infty$ then $\left\|T^{n}\right\|^{1/2} \rightarrow \infty$ and we have
\begin{center}
$\left\|T^{n}-A\right\| \geq \left\|T^{n}\right\| - \left\|A\right\| \geq \frac{1}{2}\left\|T^{n}\right\| $   \; for all n large enough.
\end{center}
 
 \begin{center}
So,\;  $\sum^{\infty}_{n=1} \frac{1}{\left\|T^{n} - A\right\|^{1/2}} < \infty$
 \end{center}
Therefore, the operators $ S_{n} :=T_{n} - A$ satisfy the conditions in [1 , Theorem 6] .\\
So there exists a dense set of pairs  $x\in X $ ;  $x^{\ast} \in X $ with:
\begin{center}
$\left|\left\langle  (T_{n} - A)x , x^{\ast}\right\rangle\right| > 0 $ and $\forall n \geq 1$ $\left|\left\langle (T_{n} - A)x , x^{\ast}\right\rangle\right|\rightarrow \infty $.\\
\end{center}
Thus there is a constant $C>0$ such that $inf_{n}\left|\left\langle (T_{n} - A)x , x^{\ast}\right\rangle\right|>0$ and we have a contradiction with the assumption that $\left\langle Ax , x^{\ast} \right\rangle \in W-Orb(T,x)^{-W}$.\\

The Hilbert space case can be proved similarly.\\

We finish this paper by posing the open question related to existence of operator weakly orbit
reflexive but not orbit reflexive.


\end{document}